\documentclass{article}

\usepackage{latexsym}
\usepackage{amssymb}
\usepackage{amsfonts}
\usepackage{amsmath}
\usepackage[all]{xy}

\newtheorem{thm}{\bf Theorem}[section]
\newtheorem{cor}[thm]{\bf Corollary}

\newtheorem{lem}[thm]{\bf Lemma}
\newtheorem{prop}[thm]{\bf Proposition}
\newtheorem{defn}[thm]{\bf Definition}

\newtheorem{rem}[thm]{\bf Remark}

\def\proof{{\parindent0pt {\bf Proof.\ }}}

\def\N{{\mathbb{N}}}

\def\pd{{\rm pd}}
\def\fd{{\rm fd}}
\def\id{{\rm id}}

\def\Did{{\rm Did}}
\def\Dpd{{\rm Dpd}}
\def\Gpd{{\rm Gpd}}
\def\PGFpd{{\rm {PGF}\mbox{-}pd}}

\def\Gfd{{\rm Gfd}}
\def\Gid{{\rm Gid}}

\def\PGFD{{\rm PGFD}}
\def\GDID{{\rm GDID}}

\def\GPD{{\rm GPD}}
\def\GID{{\rm GID}}

\def\gldim{{\rm gldim}}

\def\DPD{{\rm DPD}}
\def\DID{{\rm DID}}

\def\proj{\mathcal{P}}
\def\flat{\mathcal{F}}
\def\inj{\mathcal{I}}
\def\FPinj{\mathcal{FI}}
\def\cot{\mathcal{C}}
\def\X{\mathcal{X}}
\def\Y{\mathcal{Y}}
\def\A{\mathcal{A}}
\def\E{\mathcal{E}}

\def\M{\mathcal{M}}

\def\Q{\mathcal{Q}}
\def\W{\mathcal{W}}
\def\R{\mathcal{R}}

\def\XX{\textbf{X}}

\def\PP{\textbf{P}}

\def\PGF{\mathcal{PGF}}

\def\GF{\mathcal{GF}}
\def\GC{\mathcal{GC}}
\def\GI{\mathcal{GI}}
\def\GP{\mathcal{GP}}

\def\DP{\mathcal{DP}}

\def\GC{\mathcal{GC}}
\def\GI{\mathcal{GI}}
\def\GP{\mathcal{GP}}

\def\resdim{{\rm resdim}}
\def\coresdim{{\rm coresdim}}

\def\Coker{{\rm Coker}}
\def\Ker{{\rm Ker}}

\def\Ho{{\rm Ho}}

\def\Ext{{\rm Ext}}
\def\Tor{{\rm Tor}}
\def\Hom{{\rm Hom}}

\def\sup{{\rm sup}}

\newcommand{\cqfd}
{\rule{2mm}{2mm}%
\medbreak%
\par%
}

\begin{document}
	\baselineskip=15pt
	
	\title{Model structures, $n$-Gorenstein flat modules and PGF dimensions}
	
	\author{Rachid El Maaouy}
	
	\date{}
	
	\maketitle
	
	{\small CeReMaR Research Center, Faculty of Sciences, B.P. 1014, Mohammed V University in Rabat, Rabat, Morocco.
	
	E-mail address: rachid\_elmaaouy@um5.ac.ma; elmaaouy.rachid@gmail.com}
	
	\begin{abstract}
		Given a non-negative integer $n$ and a ring $R$ with identity, we construct a hereditary abelian model structure on the category of left $R$-modules where the class of cofibrant objects coincides with $\GF_n(R)$ the class of left $R$-modules with Gorenstein flat dimension at most $n$, the class of fibrant objects coincides with  $\flat_n(R)^\perp$ the right $\Ext$-orthogonal class of left $R$-modules with flat dimension at most $n$, and the class of trivial objects coincides with $\PGF(R)^\perp$  the right $\Ext$-orthogonal class of PGF left $R$-modules recently introduced by \v{S}aroch and \v{S}\v{t}ov\'{\i}\v{c}ek. The homotopy category of this model structure is triangulated equivalent to the stable category $\underline{\GF(R)\cap\cot(R)}$ modulo flat-cotorsion modules and it is compactly generated when $R$ has finite global Gorenstein projective dimension.
		
		The second part of this paper deals with the PGF dimension of modules and rings. Our results suggest that this dimension could serve as an alternative definition of the  Gorenstein projective dimension. We show, among other things, that ($n$-)perfect rings can be characterized in terms of Gorenstein homological dimensions, similar to the classical ones, and the global Gorenstein projective dimension coincides with the global PGF dimension. 
		
	\end{abstract}
	\medskip
	
	
	
	\normalsize
	\section{Introduction}

	\hskip .5cm Throughout this paper, $R$ will be an associative ring with identity, and all modules will be, unless otherwise specified, unital left $R$-modules. When right $R$-modules need to be used, they will be denoted as $M_R$, while in these cases left $R$-modules will be denoted by $_R M$. whenever it is convenient, right $R$-modules are identified with modules over the opposite ring $R^{op}$.  We use $\inj(R)$, $\proj(R)$, $\flat(R)$ and $\cot(R)$ to denote the classes of injective, projective,  flat and cotorsion $R$-modules, respectively.

	From a theorem by Christensen, Estrada and Thompson \cite[Theorem 4.5]{CET20}, the subcategory $\GF(R)\cap\cot(R)$ of Gorenstein flat and cotorsion $R$-modules is a Frobenius category with projective-injective objects all flat-cotorsion $R$-modules. Hence, the stable category $\mathcal{E}_0:=\underline{\GF(R)\cap\cot(R)}$ is a triangulated category. On the other hand, \v{S}aroch and  \v{S}\v{t}ov\'{\i}\v{c}ek  \cite[Section 4]{SS20} constructed a hereditary abelian model structure $\M_0$ where the class of cofibrant objects coincides with $\GF(R)$ and the class of fibrant objects coincides with $\cot(R)$. As a particular case of Estrada, Iacob, and  P\'{e}rez \cite[Corollary 4.6]{EIP20} (see also \cite[Corollary 4.18]{BEGO22b}) the homotopy category $\Ho(\M_0)$ of this model structure is triangulated equivalent to the stable category  $\E_0$.
	
	In the first part of this paper, Section 3, we are interested in a general situation of the same kind of problem: that of modules of bounded Gorenstein flat dimension. It is therefore natural to ask whether we can construct a model structure such that the class of cofibrant objects coincides with the class $\GF_n(R)$ of all modules having Gorenstein flat dimension at most an integer $n\geq 0$.
	
	First, we show that there is a trivial model structure (in the sense that every module is a trivial object) $$\mathcal{N}_n=(\GF_n(R), R\text{-Mod},\GF_n(R)^\perp).$$
	In fact, this is nothing but saying that $\GF_n(R)$  is the left hand of a complete cotorsion pair. However, the following result shows more than that.
	
	\bigskip
	\noindent{ \bf Theorem A.} \textbf{(The $n$-Gorenstein flat cotorsion pair)} The class $\GF_n(R)$ forms the left hand of a perfect hereditary cotorsion pair $(\GF_n(R),\GC_n(R))$ with $\GC_n(R)=\PGF(R)^\perp\cap \flat_n(R)^\perp$.
	
	In particular, $\GF_n(R)$ is covering and $\GC_n(R)$ is enveloping.
	
	\bigskip
	
	Theorem A is essential to obtain our next main result (Theorem \ref{GFn mod struc}). This result provides us with a (non-trivial) model structure  where the class of fibrant objects coincides with $\cot_n(R):=\flat_n(R)^\perp$ the right $\Ext$-orthogonal class of all $R$-modules of flat dimension at most $n$.

	One of the consequences of this model structure is that the stable category $\E_n=\underline{\GF_n(R)\cap\cot_n(R)}$ modulo objects in $\flat_n(R)\cap\cot_n(R)$ is triangulated (Corollary \ref{quas-Frob}). In the case $n=0$,  Liang  and Wang \cite[Theorem 1(a)]{LW20} have recently shown that this stable category is compactly generated when $R$ is right coherent with finite global Gorenstein projective dimension. It is natural to ask whether this generatedness property is inherited by $\mathcal{E}_n$ for all $n\geq 0$. Based on their work, we answer this question in the positive, assuming only that $R$ has finite global Gorenstein AC-projective dimension.
	
	\bigskip
	
	\noindent{ \bf Theorem B.} (\textbf{The $n$-Gorenstein flat model structure}) \label{ThA} There exists a hereditary abelian model structure on $R$-Mod $$\M_n=\left(\GF_n(R),\PGF(R)^\perp, \cot_n(R)\right).$$
	
	Consequently, there are triangulated equivalences
	$$\underline{\GF_n(R)\cap \cot_n(R)}\simeq\cdots \simeq \underline{\GF(R)\cap \cot(R)}\simeq \underline{\PGF(R)}.$$
	
	Furthermore, if $R$ has finite global Gorenstein  AC-projective dimension, then these triangulated categories are compactly generated.
	
	\bigskip
	
	The $n$-Gorenstein flat model structure has been found for two particular rings: (Iwanaga-)Gorenstein rings by P\'erez \cite[Theorem 14.3.1]{Per16} where the trivial objects are known and coincide with  modules having finite flat dimension and right coherent rings by A. Xu \cite[Theorem 4.5(2)]{Xu17} where less information is known about the class of trivial objects.
	
	Therefore, Theorem B generalizes and improves their results in the sense that we do not need any assumption on the ring $R$ with an explicit description of the trivial objects. This makes Theorem B of interest since the class of trivial objects is the most important class in any model structure for it determines the corresponding homotopy category, as explained in the fundamental theorem of model categories \cite[Theorem 1.2.10]{Hov99}. 
	
	We point out that the proof of Theorem B is different from that of P\'erez and Xu. It is mainly based on the work developed by \v{S}aroch and  \v{S}\v{t}ov\'{\i}\v{c}ek in \cite{SS20}. One of the key results (Theorem \ref{Charac of GFn}) is a new and useful characterization  of modules having finite Gorenstein flat dimension.

	$$*\;*\;*$$
	
	A guiding principle in the development of Gorenstein homological algebra is to study analogues of results concerning absolute homological dimensions.  For example, consider the well-known questions:  (Q1) Is every Gorenstein projective module Gorenstein flat? and (Q2) Is the class of Gorenstein projective modules special precovering?  These questions have been studied by many authors, such as  Cort\'{e}s-Izurdiaga and \v{S}aroch \cite{CS21}, Emmanouil \cite{Emm12}, Enochs and Jenda \cite{EJ00} and  Iacob \cite{Iac20} among others. One could also ask similar questions in the case of Ding projective modules. However, all these questions remain open.

	\v{S}aroch and \v{S}\v{t}ov\'{\i}\v{c}ek \cite{SS20}, on the other hand, introduced PGF modules and showed that they are Gorenstein flat and they form a special precovering class. These properties can be seen as positive answers to Questions (Q1) and (Q2) if we think of PGF modules as an alternative definition of Gorenstein projective modules.

	Our main purpose of Section 4 is to support the following claim: \textbf{"the PGF dimension could serve as an alternative definition of the Gorenstein (Ding) projective dimension over any ring"}.

	Recall \cite[Definition 2.1]{EJL05} that a ring $R$ is left $n$-perfect if every flat $R$-module has projective dimension $\leq n$. In particular, left perfect rings are exactly left $0$-perfect rings. We show (Theorem \ref{n-perf}) that a ring $R$ is left $n$-perfect if and only if every Gorenstein flat left $R$-module has PGF dimension at most $n$. In the case $n=0$, we obtain other equivalent assertions similar to the classical ones, supporting the above claim.

	It follows by \cite[Theorem 4.4]{SS20} that any PGF module is Gorenstein projective. But whether the converse is true remains open. This leaves us with another question: (Q3) When is any Gorenstein projective a PGF module? This question was first asked and investigated by \v{S}aroch and  \v{S}\v{t}ov\'{\i}\v{c}ek \cite{SS20} and later by Iacob \cite{Iac20}. It turns out that this question is closely related to Questions (Q1)  and (Q2) (see Remark \ref{Q1 and Q2}).
	
	One nice way to study Question (Q3) is to measure how far a Gorenstein projective module is from being PGF. It is proved (Proposition \ref{PGF=GP}) that this can be reduced to two cases: either Gorenstein projective modules are as close as possible to PGF, i.e., Gorenstein projective and PGF modules coincide, or, as far as possible, and this means that the PGF dimension of any Gorenstein projective module is infinite.

	Next, we obtain a variety of conditions that are equivalent to the first case (Proposition \ref{PGF=GP}). Under these equivalent assertions, we get a positive answer to Questions (Q1) and (Q2) (Corollary \ref{Q1 and Q2}). In particular, this is the case when $R$ is right weak coherent and left $n$-perfect or $R$ is a ring such that every injective right $R$-module has finite flat dimension.

	Right weak coherent rings were introduced by Izurdiaga in \cite{CI16} as a natural generalization of right coherent rings. They are rings for which the direct product of any family of flat $R$-modules has finite flat dimension. It turns out that these rings provide a general framework in which one can replace the assumption of being coherent with that of being weak coherent (see Corollary \ref{PGF=DP} and Remark \ref{coh to weak coh}). 
	
	The rest of Section 4 is devoted to the global dimension of $R$ with respect to the class $\PGF(R)$ and its link to other global dimensions. First, we provide simple ways to compute it (Theorem \ref{charc of PGFD}). This result is then used to show our last main result of this paper which, along with its consequences (Corollaries \ref{PGFD=GID=DID} and \ref{Aus}), clearly support our claim above. The following result states that, unlike the classes $\PGF(R)$ and $\GP(R)$, their global dimensions coincide. 
	
	\bigskip 
	\noindent{ \bf Theorem C.} For any ring $R$, we have the following equality:  $$\sup\{\PGFpd_R(M)|\text{ $M$ is an $R$-module}\}=\sup\{\Gpd_R(M)|\text{ $M$ is an $R$-module}\}.$$

	\section{Preliminaries}

	\noindent{\bf Resolutions.} Given a class  $\X$ of $R$-modules and a class $\Y$ of right $R$-modules, an $\X$-resolution of an $R$-module $M$ is an exact complex $\cdots\to X_1\to X_0\to M\to 0$ where
	$X_i\in\X$.

	A sequence $\XX$ of $R$-modules is called $\left(\Y\otimes_R-\right)$-exact (resp., $\Hom_R(\X,-)$-exact,
	$ \Hom_R(-,\X)$-exact) if $Y\otimes_R\XX$  (resp., $\Hom_R(X,\XX)$, $\Hom_R(\XX,X)$)
	is an exact complex for every $Y\in\Y$ (resp., $X\in\X$).

	An $R$-module $M$ is said to have $\X$-resolution dimension at most an integer $n\geq 0$, $\resdim_\X(M)\leq n$, if $M$ has a finite $\X$-resolution: $0\to X_n\to\cdots\to X_1\to X_0\to M\to 0.$ If $n$ is the least non negative integer for which such a sequence exists then its $\X$-resolution dimension is precisely $n,$ and if there is no such $n$ then we define its $\X$-resolution dimension as $\infty$.

	Given a class of $R$-modules $\mathcal{Z}$, the $\X$-dimension of $\mathcal{Z}$ is defined as: 
	$$\resdim_\X(\mathcal{Z})=\sup\{\resdim_\X(Z)|Z\in\mathcal{Z}\}.$$
	$\X$-coresolutions and  $\X$-coresolution dimensions are defined dually.
	
	\noindent{\bf Gorenstein modules and dimensions.} An $R$-module $M$ is called  Gorenstein flat if  it is a syzygy of an  $(\inj(R^{op})\otimes_R-)$-exact exact sequence of flat $R$-modules.  Replacing flat with projective, we get the definition of projectively coresolved Gorenstein flat (PGF for short) modules  \cite{SS20}. An $R$-module $M$ is called Gorenstein projective if it is a syzygy of a $\Hom_R(-,\proj(R))$-exact exact sequence of projective $R$-modules. Gorenstein injective modules are defined dually. We let $\GF(R)$
	(resp., $\PGF(R)$, $\GP(R)$, and $\GI(R)$) denote the subcategory of Gorenstein flat (resp., PGF, Gorenstein projective, and Gorenstein injective) $R$-modules.
	
	The Gorenstein flat (resp., PGF, Gorenstein projective, Gorenstein injective) dimension of an $R$-module $M$ is defined as  $\Gfd_R(M):=\resdim_{\GF(R)}(M)$ (resp., $\PGFpd_R(M):=\resdim_{\PGF(R)}(M)$, $\Gpd_R(M):=\resdim_{\GP(R)}(M)$, and $\Gid_R(M):=\coresdim_{\GI(R)}(M)$).

	Let $\GPD(R)$ (resp, $\GID(R)$) denote the global Gorenstein projective (resp., global Gorenstein injective) dimension of $R$. i.e., $$\GPD(R):=\Gpd_R(R\text{-Mod})=\sup\{\Gpd_R(M)|M\in R\text{-Mod}\}$$ (resp., $\GID(R):=\Gid_R(R\text{-Mod})=\sup\{\Gid_R(M)|M\in R\text{-Mod}\}$).

	\begin{defn} Let $\flat_n(R)$ and $\GF_n(R)$ denote the subcategories  of $R$-modules having flat and Gorenstein flat dimension at most $n\geq 0$. We call modules in  these subcategories  $n$-flat and  $n$-Gorenstein flat, respectively.

    An $R$-module $M$ is said to be $n$-cotorsion (resp., $n$-Gorenstein cotorsion) if $\Ext^1_R(N,M)$ for all $n$-flat (resp., $n$-Gorenstein flat) modules $N$. We use  $\cot_n(R)$ and $\GC_n(R)$ to denote the subcategrories of $n$-cotorsion and $n$-Gorenstein cotorsion modules, respectively.
		
	\end{defn}
	Note that $0$-flat (resp., Gorenstein $0$-flat, $0$-cotorsion, $0$-Gorenstein cotorsion) modules coincide with the flat (resp., Gorenstein flat, cotorsion, Gorenstein cotorsion) modules.   
	
	\bigskip
	
	\noindent{\bf Cotorsion pairs.} Given an abelian category $\A$, and a class of objects $\X$ in $\A$, we use the following standard notations: $\X^\perp=\{A\in\A|\Ext^1_\A(X,A)=0,\forall X\in\X\}$ and 
	$^\perp\X=\{A\in\A|\Ext^1_\A(A,X)=0,\forall X\in\X\}.$
	
	An $\X$-precover of an object $A\in\A$ is a morphism $f:X\to A$ with $X\in \X$, in such a way that $f_*:\Hom_\A(X',X)\to \Hom_\A(X',M)$ is surjective for every $X'\in \X$. An $\X$-precover is called an $\X$-cover if every endomorphism $g: X\to X$ such that $fg=f$ is an automorphism of $X$. If every object has an $\X$-(pre)cover then the class $\X$ is said to be  (pre)covering.  An $\X$-precover is called special if it is epimorphism and $\Ker f\in\X^\perp$. $\mathcal{X}$-(special pre)envelopes can be defined dually.

	A pair $\left(\X,\Y\right)$ of classes of objects in $\A$ is called a cotorsion pair if $\X^\perp=\Y$ and $\X=^\perp\Y$. 
	A cotorsion pair $\left(\X,\Y\right)$ is said to be: 
	
	$\bullet$ Hereditary if $\Ext^i_\A(X,Y)=0$ for every $X\in\X$, $Y\in\Y$ and  $i\geq 1$. 
	
	$\bullet$ Complete if any object in $\A$ has a special $\X$-precover  and a special $\Y$-preenvelope. 
	
	$\bullet$ Perfect if every module has an $\X$-cover and a $\Y$-envelope. 
	
	$\bullet$ Cogenerated by a set if there is a set $\mathcal{S}$ such that $\X=\mathcal{S}^\perp$.

	It is well known that a perfect cotorsion  pair $(\X,\Y)$ in $R$-Mod is complete. The converse holds when $\X$ is closed under direct limits \cite[Corollary 2.3.7]{GT12}.  A well-known  method for constructing complete cotorsion pairs in $R$-Mod is to cogenerate one from a set (see \cite[Theorem 3.2.1]{GT12}, for instance). 
	
	\bigskip
	\noindent{\bf Abelian model structures.} 
	Given a bicomplete category $\A$, a model structure on $\A$ is given by three classes of morphisms of $\A$, called cofibrations, fibrations and weak equivalences that satisfy a set of axioms \cite[Definition 1.1.3]{Hov99}. 
	
	Over a (bicomplete) abelian category $\A$, Hovey defined abelian model structures and showed that we can shift all our attention from morphisms to objects \cite[Theorem 2.2]{Hov02}. Namely, an abelian model structure on $\A$, is equivalent to a triple (known as Hovey triple) $\M=(\Q,\W,\R)$ of classes of objects in $\A$ such that $(\Q, \W\cap\R)$ and $(\Q\cap\W,\R)$ are complete cotorsion pairs and $\W$ is thick.  In this case, $Q$ (resp., $\R$, $\W$) is precisely the class of cofibrant (resp., fibrant, trivial) objects and an abelian model structure (in the sense of \cite[Definition 1.1.3]{Hov99}) is one which the following two conditions are satisfied: 
	
	(i) A morphism is a (trivial) cofibration if and only if it is a monomorphism with (trivially) cofibrant cokernel.
	
	(ii) A morphism is a (trivial) fibration if and only if it is an epimorphism with (trivially) fibrant kernel.
	
	An abelian model structure is called hereditary if both of the associated cotorsion pairs are hereditary (see Gillespie \cite{Gil16} for a nice survey on hereditary abelian model categories).  In this paper, we often identify abelian model structures  with Hovey triples. 
	
   Given a model structure $\M$, the homotopy category, denoted as $\Ho_\A(\M)$, is defined by formally inverting the weak equivalences of $\M$. More precisely, $\Ho_\A(\M)$ is obtained after localizing $\A$ at the class of weak equivalences. If moreover, $\M$ is hereditary, the homotopy category is known to be an algebraic triangulated category  and it encodes a variety of (relative) homological algebra on $\A$.

	\section{n-Gorenstein flat modules and their model structure}
	
	\hskip .5cm In this section, we aim at proving Theorems A and B in the introduction. First, we need some basic properties of $n$-Gorenstein flat modules.
	
	\begin{prop}\label{Gn-flat prop} The class of $n$-Gorenstein flat $R$-modules is projectively resolving, closed under direct sums, direct summands, and direct limits.
	\end{prop}
	\proof  It is clear that any projective module is Gorenstein flat, and then $n$-Gorenstein flat. On the other hand, taking $W=R$ in \cite{BEGO22}, we get that $\GF(R)={\rm G_WF}(R)$ and this class is closed under extensions by \cite[Theorem 4.11]{SS20}. Therefore, the class $\GF_n(R)$ is closed under extensions and kernels of epimorphisms by \cite[Proposition 7.9]{BEGO22} and  direct sums and summands by \cite[Corollary 7.8(2)]{BEGO22}. It remains to show that $\GF_n(R)$ is closed under direct limits. 
	
	Following  \cite[Corollary 1.7 and the Remark that follows it]{AR94} it suffices to show that the class $\GF_n(R)$ is closed under well ordered continuous direct limits of  $n$-Gorenstein flat $R$-modules: $$(M_\alpha)_{\alpha<\lambda}:= M_0\to M_1\to \cdots\to M_\omega\to M_{\omega+1}\to \cdots$$
	Here continuous (or smooth in the sense of \cite{AR94}) means that for each limit ordinal $\beta<\lambda$, $M_\beta=\varinjlim_{\alpha<\beta} M_\alpha$.
	
	$\bullet$ If $\lambda < \omega$, then $\lim_{\alpha <\lambda} M_\alpha=M_{\lambda-1}\in \GF_n(R)$.
	
	$\bullet$ Assume now that $\lambda=\omega$. By \cite[Proposition 7.11]{BEGO22}, there exists an exact sequence $0\to M_0\to F_0\to G_0\to 0$
	with $F_0\in\flat_n(R)$ and $G_0$ Gorenstein flat. Consider the push out
	$$\xymatrix{
		&0 \ar[r]  &M_0 \ar[r] \ar[d]  &F_0 \ar[r] \ar[d] & G_0 \ar[r] \ar@{=}[d] &0\\
		&0 \ar[r]  &M_1 \ar[r]   &V \ar[r]   & G_0 \ar[r]  &0.
	}$$
	Since $M_1\in\GF_n(R)$ and $G_0$ is Gorenstein flat, $V\in\GF_n(R)$. Using again \cite[Proposition 7.11]{BEGO22}, we get an exact sequence $0\to V\to F_1\to G\to 0$ with $F_1\in\flat_n(R)$ and $G$ Gorenstein flat. 
	
	Consider another push out
	
	$$\xymatrix{
		& & &0 \ar[d] &0 \ar[d]\\
		&0 \ar[r]  &M_1 \ar[r] \ar@{=}[d]  &V \ar[r] \ar[d] & G_0 \ar[r] \ar[d] &0\\
		&0 \ar[r]  &M_1 \ar[r]   &F_1 \ar[r]\ar[d]   & G_1 \ar[r]\ar[d]  &0\\
		& & &G \ar[d] \ar@{=}[r]  & G\ar[d]  \\
		& & &0   &0.}$$ 
	Since $G_0$ and $G$ are Gorenstein flat, so is $G_1$. Therefore, we get the following morphism of exact sequences induced by
	the morphism $M_0\to M_1$: 
	$$\xymatrix{
		&A_0=:\ar[d] 	&0 \ar[r]  &M_0 \ar[r] \ar[d]  &F_0 \ar[r] \ar[d] & G_0 \ar[r] \ar[d] &0\\
		&A_1=:	&0 \ar[r]  &M_1 \ar[r]   &F_1 \ar[r]   & G_1 \ar[r]  &0
	}.$$
	Repeating this process, we get a commutative diagram with exact rows
	$$\xymatrix{
		&A_0=:\ar[d]    & 0 \ar[r]  &M_0 \ar[r] \ar[d]  &F_0 \ar[r] \ar[d] & G_0 \ar[r] \ar[d] &0\\
		&A_1=:\ar[d]  &0 \ar[r]   &M_1 \ar[r]\ar[d]   &F_1 \ar[r]\ar[d]   & G_1 \ar[r]\ar[d]  &0\\
		&A_2=:\ar[d]  &0 \ar[r]  &M_2 \ar[r]\ar[d]   &F_2 \ar[r] \ar[d]  & G_2 \ar[r]\ar[d]  &0\\
		&:  &  &:   &:   & :  
	}$$
	where $F_i\in\flat_n$ and $G_i$ is Gorenstein flat for every $i\geq 0$. Applying the exact functor $\varinjlim$ to this commutative diagram, we get the induced exact sequence
	$$\varinjlim A_m=:  0\to \varinjlim M_m\to \varinjlim F_m\to \varinjlim G_m\to 0$$ 
	with $\varinjlim G_m\in \GF(R)\subseteq \GF_n(R)$ by \cite[ Lemma 3.1]{YL12} and  $\varinjlim F_m\in \flat_n(R)$ since the functor $\Tor^R_{n+1}(X,-)$ commutes with direct limits for any right $R$-module $X$. Hence, $M_\omega=\varinjlim M_m\in\GF_n(R)$.
	
	Finally, using transfinite induction, it is clear that the above argument generalizes and we get that $\varinjlim_{\alpha<\lambda} M_\alpha\in\GF_n(R)$ as desired.
	\cqfd

	The following two lemmas will be used in several places in this paper. 
	\begin{lem} \label{n-flat is PGC} The inclusion
		$\flat_n(R) \subseteq   \PGF(R)^\perp$ holds true. 
		
	\end{lem}
	\proof  It follows from the fact that $ \PGF(R)^\perp$ is thick by \cite[Theorem 4.9]{SS20} and  $\flat(R)\subseteq \PGF(R)^\perp$ by \cite[Theorem 4.11]{SS20}. \cqfd

	\begin{lem} \label{n-flat cot-pair} $(\flat_n(R),\cot_n(R))$ is a perfect and hereditary cotorsion pair cogenerated by a set.
	\end{lem}
	\proof Follows by \cite[Theorem 3.4(2)]{DM07} and its proof \cqfd
	
	The following characterization of $n$-Gorenstein flat modules is a key result to show many results in this section. The case $n=0$ is due to  \v{S}aroch and \v{S}\v{t}ov\'{\i}\v{c}ek \cite[Theorem 4.11]{SS20}. 
	
	\begin{thm} \label{Charac of GFn} The following assertions are equivalent for any $R$-module $M$.
		\begin{enumerate}
			\item $M$ is $n$-Gorenstein flat.
			
			\item There exists a $\Hom_R(-,\flat_n(R)\cap\cot_n(R))$-exact exact sequence $$0\to K\to L\to M\to 0$$
			where $K\in \flat_n(R)$ and $L\in \PGF(R)$.
			
			\item $\Ext^1_R(M,C)=0$ for every $C\in  \PGF(R)^\perp\cap \cot_n(R)$.
			
			\item There exists an exact sequence of $R$-modules $$0\to M\to F\to N\to 0$$
			where $F\in\flat_n(R)$ and $N\in \PGF(R)$.
			
		\end{enumerate}
		
		Consequently, $\flat_n(R)=\GF_n(R)\cap \PGF(R)^\perp$. 
	\end{thm}

	\proof  $4.\Rightarrow 1.$ It follows by Proposition \ref{Gn-flat prop} since $F\in\flat_n(R)\subseteq \GF_n(R)$ and $N\in \PGF(R)\subseteq \GF_n(R)$.
	
	$1.\Rightarrow 4.$ Assume that $\Gfd_R(M)\leq n$ and proceed by induction on n. 
	
	The case $n=0$ follows from \cite[Theorem 4.11(4)]{SS20}.
	
	Assume that  $n\geq 1$. There exists an exact sequence  $0\rightarrow G_n\rightarrow\cdots\rightarrow G_0\rightarrow M\rightarrow 0$, where $G_i\in \GF(R)$ for every $i=0,\cdots,n.$ Let $K=\Ker(G_0\to M)$. Clearly, $\Gfd_R(K)\leq n-1$. So, by induction there exists an exact sequence $0\rightarrow K\rightarrow F'\rightarrow N'\rightarrow 0,$ where $N'\in \PGF(R)$ and $\fd_R(F')\leq n-1$. Consider the pushout diagram $$\xymatrix{ & 0\ar[d] & 0\ar[d] & & \\  0\ar[r] & K\ar[d] \ar[r] & G_0\ar@{-->}[d] \ar[r] & M \ar@{=}[d] \ar[r] & 0 \\ 0\ar[r] & F'\ar@{-->}[r]\ar[d] & D \ar[d] \ar[r] & M\ar[r] & 0 \\ & N'\ar@{=}[r] \ar[d] & N'\ar[d] & & \\ & 0 & 0}.$$
	By the middle column $D$ is Gorenstein flat. Using induction again, we get a short exact sequence of $R$-modules  $0\to D\to F_0\to N\to 0$ where $F_0$ is flat and $N\in\PGF(R)$. 
	
	Consider now another pushout diagram 
	
	$$\xymatrix{ & & 0\ar[d] & 0\ar[d] & \\ 0\ar[r] & F'\ar@{=}[d] \ar[r] & D\ar[d] \ar[r] & M\ar@{-->}[d] \ar[r] & 0 \\ 0\ar[r] & F'\ar[r] & F_0 \ar@{-->}[r]\ar[d] & F\ar[r] \ar[d] & 0 \\ & & N\ar[d]\ar@{=}[r] & N\ar[d] & \\ & & 0 & 0}$$
	It remains to see that $\fd_R(F)\leq n$. But this is true by the middle row since $F_0$ is flat and $\fd_R(F')\leq n-1$. 
	
	$1.\Rightarrow 2.$  The case $n=0$ holds by  \cite[Theorem 4.11(2)]{SS20}. So, we may assume that $n\geq 1$. By \cite[Theorem 4.9]{SS20}, there exists a special PGF precover $0\to N\to G\to M\to 0.$ By \cite[Proposition 7.9]{BEGO22}, $\Gfd_R(N)\leq \sup\{\Gfd_R(G),\Gfd_R(M)-1\}\leq n-1$ which implies by  $(4)$ that there exists an exact sequence of $R$-modules $0\to N\to K\to G'\to 0$
	where $\fd_R(K)\leq n-1$ and $G'\in\PGF(R)$. Consider the pushout
	$$\xymatrix{ & 0 \ar[d] & 0 \ar[d] \\ 0 \ar[r] & N \ar[d] \ar[r] & G \ar[r] \ar@{-->}[d] & M \ar[r] \ar@{=}[d] & 0 \\ 0 \ar[r] & K \ar@{-->}[r] \ar[d] & L \ar[d] \ar[r] & M \ar[r] & 0 \\ & G' \ar@{=}[r]\ar[d] & G' \ar[d] \\ & 0 & 0}$$
	Since $\PGF(R)$ is closed under extensions, $L\in\PGF(R)$. Then, there exits (by definition) an exact sequence $0\to L\to P\to L'\to 0$,  with $P$ projective and $L'\in \PGF(R)$. Consider now another pushout
	$$\xymatrix{ &  & 0 \ar[d] & 0 \ar[d] \\ 0 \ar[r] & K \ar@{=}[d] \ar[r] & L \ar[r] \ar[d] & M \ar[r] \ar@{-->}[d] & 0 \\ 0 \ar[r] & K \ar[r]  & P \ar[d] \ar@{-->}[r] & H \ar[r]\ar[d] & 0 \\ & & L' \ar[d]\ar@{=}[r]  & L' \ar[d]\\ & & 0 & 0}.$$
	Since $\fd_R(K)\leq n-1$ and $P$ is projective (then flat),  $H\in\flat_n(R)$ by the middle row. It follows from the following commutative digram with exact rows 
	$$\xymatrixcolsep{1.7pc}\xymatrix{ 0 \ar[r] &\Hom_R(H,X) \ar[r] \ar[d] &\Hom_R(P,X) \ar[r] \ar[d]  &\Hom_R(K,X) \ar@{=}[d] \ar[r] & \Ext^1_R(H,X)=0\\
		0\ar[r] &\Hom_R(M,X) \ar[r] 	&\Hom_R(L,X) \ar[r]     & \Hom_R(K,X)
	}$$
	that $0\to L\to K\to M\to  0$ is $\Hom_R(-,X)$-exact for any module $X\in\flat_n(R)\cap\cot_n(R)$. 
	
	$ 2.\Rightarrow 3.$  Similar to \cite[Theorem 4.11(2)$\Rightarrow$(3)]{SS20}, using the fact that any module $C\in \PGF(R)^\perp\cap \cot_n(R)$  has a special $n$-flat precover which exists by  Lemma \ref{n-flat cot-pair}.

	$ 3.\Rightarrow 4.$ Similar to \cite[Theorem 4.11(3)$\Rightarrow$(4)]{SS20} using the fact that every module $F\in \PGF(R)^\perp$ has a special $n$-flat precover which exists by Lemma \ref{n-flat cot-pair} and that $\flat_n(R)\subseteq \PGF(R)^\perp$ by Lemma \ref{n-flat is PGC}. 
	
	Finally, it remains to show the equality. The inclusion  $\flat_n(R)\subseteq \GF_n(R)$ is clear and by Lemma \ref{n-flat is PGC} we also have $\flat_n(R)\subseteq \PGF(R)^\perp$. Then, $\flat_n(R)\subseteq \GF_n(R)\cap \PGF(R)^{\perp}$. Conversely, if $M\in\GF_n(R)\cap \PGF(R)^{\perp}$, then there exists by (4) a split exact sequence of $R$-modules  $0\to M\to F\to N\to 0$, where $N\in\PGF(R)$ and $\fd_R(X)\leq n$. Then, $M\in \flat_n(R)$. \cqfd

	The following result, which is Theorem A from the introduction, improves and  generalizes \cite[Proposition 14.3.5]{Per16} and \cite[Lemma 4.4]{Xu17}. Here we have a  more explicit description of the right hand of the $n$-Gorenstein flat cotorsion pair without any assumptions on the ring.
	\begin{cor}\label{GFn cot pair} The pair $\left(\GF_n(R),\PGF(R)^\perp\cap \cot_n(R)\right)$ is a perfect hereditary cotorsion pair cogenerated by a set. 
		
		Consequently, the following assertions hold.
		\begin{enumerate}
			\item[(a)] Every $R$-module has a $n$-Gorenstein flat cover and $n$-Gorenstein cotorsion envelope. 
			\item[(b)] An $R$-module is $n$-Gorenstein cotorsion if and only if it is $n$-cotorsion from $\PGF(R)^\perp$.
		\end{enumerate}
	\end{cor}
	\proof Let $\mathcal{A}=\GF_n(R)$ and $\mathcal{B}=\PGF(R)^\perp \cap \cot_n(R)$. By Theorem \ref{Charac of GFn}(3),  $\mathcal{A}=\;^{\perp}\mathcal{B}$. On the other hand, clearly $\PGF(R)\cup \flat_n(R)\subseteq \GF_n(R)$. Then, $\mathcal{A}^{\perp}=\GF_n(R)^{\perp}\subseteq \flat_n(R)^\perp \cap \PGF(R)^\perp=\mathcal{B}$. Hence, $\mathcal{A}^{\perp}=\mathcal{B}$ and $(\mathcal{A},\mathcal{B})$ is a cotorsion pair.
	
	Now we prove that this cotorsion pair is perfect. But, since $\GF_n(R)$ is closed under direct limits by Proposition \ref{Gn-flat prop}(3), we only need to check that it is complete \cite[Corollary 2.3.7]{GT12}. Following \cite[Theorem 3.2.1(b)]{GT12}, a nice way to do that is to show that it is cogenerated by a set. By the proof of \cite[Theorem 4.9]{SS20}, $(\PGF(R),\PGF(R)^\perp)$ is a cotorsion pair cogenerated by a set  $\mathcal{A}_1$. We also know by Lemma \ref{n-flat cot-pair} that the cotorsion pair $(\flat_n(R),\cot_n(R))$ is cogenerated by a set $\mathcal{A}_2$. Then, $\PGF(R)^\perp \cap \cot_n(R)=\mathcal{A}_1^\perp \cap \mathcal{A}_2^\perp=(\mathcal{A}_1\cup\mathcal{A}_2)^\perp$. This means that our desired pair is cogenerated by a set and hence complete. \cqfd

	Now we are ready to prove our second main result in this section. The following result improves \cite[Theorem 14.3.1]{Per16} and  \cite[Theorem 4.5]{Xu17} in the sense that we need no extra assumptions on the ring $R$ with an obvious description of the trivial objects. 
	\begin{thm}\label{GFn mod struc} For any ring $R$, there exists a hereditary abelian model structure  on $R$-Mod $$\mathcal{M}_n:=\left(\GF_n(R),\PGF(R)^\perp, \cot_n(R)\right).$$ Equivalently, there exists a hereditary abelian model structure on $R$-Mod, in which the  cofibrations coincide with the monomorphisms with $n$-Gorenstein flat cokernel, the fibrations coincide with the epimorphisms with $n$-cotorsion  kernel, and $\PGF(R)^\perp$ is the class of trivial objects.

	\end{thm}
	\proof By Theorem \ref{Charac of GFn}, $\flat_n(R)=\GF_n(R)\cap \PGF(R)^\perp$. Then,  by Lemma \ref{n-flat cot-pair} and Corollary \ref{GFn cot pair}, there are two complete and hereditary cotorsion pairs:
	$$(\GF_n(R),\PGF(R)^\perp\cap\cot_n(R)) \text{ and } (\GF_n(R)\cap \PGF(R)^\perp,\cot_n(R)).$$
	Finally, the class $\PGF(R)^\perp$ is thick by \cite[Theorem 4.9]{SS20} which completes the proof. \cqfd
	
	Recall that a Frobenius category $\A$ is an exact category with enough injectives and projectives such that the projective objects coincide with the injective objects. In this case, we can form the stable category $\underline{\A}:=\A/\sim$, which has the same objects as $\A$ and $\Hom_{\A/\sim}(X,Y):=\Hom_{\A}(X,Y)/\sim$, where $f\sim g$ if and only if $f-g$ factors through a projective-injective object.
	
	The following result can be proved directly as in the case $n=0$ (see \cite[Theorem 4.5]{CET20}). However, we follow an alternative approach due to Gillespie \cite{Gil11}. Recall \cite[Sections 4 and 5]{Gil11}  that for any hereditary Hovey triple $\M=(\Q,\W,\R)$, the subcategory $\Q\cap \R$ is a Frobenius exact category (with the exact structure given by the short exact sequences with terms in $\Q\cap \R$)  whose projective-injective objects are precisely those in $\Q\cap\W\cap\R$. Moreover, the stable category $\underline{\Q\cap\R}$ is triangulated equivalent to $\Ho(\M)$ (see also \cite[Proposition 4.2 and Theorem 4.3]{Gil16}). Consequently, we have the following corollary.
	\begin{cor}\label{quas-Frob} The subcategory  $\GF_n(R)\cap \cot_n(R)$ of both $n$-Gorenstein flat and $n$-cotorsion $R$-modules is a Frobenius category. The projective-injective modules are given by objects in $\flat_n(R)\cap\cot_n(R)$. Moreover, the homotopy category of the $n$-Gorenstein flat model structure is triangle equivalent to the stable category
		$\underline{\GF_n(R)\cap \cot_n(R)}$.
	\end{cor}
	
	Let $(\mathcal{T},\Sigma)$ be a triangulated category with coproducts and an autofunctor $\Sigma$. Recall that an object $C$ of $\mathcal{T}$ is compact if for each family $\{Y_i|i\in I\}$ of objects of $\mathcal{T}$, the canonical morphism
	$$\bigoplus_{i\in I} \Hom_\mathcal{T}(C,Y_i)\to  \Hom_\mathcal{T}(C,\bigoplus_{i\in I}Y_i)$$
	is an isomorphism. The category $\mathcal{T}$ is called compactly generated if there exists a set $\mathcal{S}\subseteq \mathcal{T}$ of 
	compact objects such that for each $0\neq Y\in\mathcal{Y}$ there is a morphism $0\neq f:\Sigma^m S\to Y$ for some $S\in\mathcal{S}$ and $m\in\mathbb{Z}$.

	It follows by Liang  and Wang \cite[Theorem 1(a)]{LW20} that the stable category $\E_0=\underline{\GF(R)\cap\cot(R)}$ is compactly generated when $R$ is right coherent with finite left global Gorenstein projective dimension. Based on their work, this property is inherited by the stable categories $\underline{\PGF(R)}$ and $\underline{\GF_n(R)\cap\cot_n(R)}$ for every integer $n\geq 0$ under the weaker assumption that $R$ has finite global Gorenstein AC-projective dimension. That is, the supremum of the  $\GP_{AC}(R)$-resolution dimension of all $R$-modules, 
     $$\GPD_{AC}(R):=\sup\{ \resdim_{\GP_{AC}(R)}(M)| \text{ $M$ is an $R$-module} \},$$
     is finite, where $\GP_{AC}(R)$ denotes the class of Gorenstein AC-projective modules in the sense of \cite{BGH14}.

	\begin{cor} \label{gener}There exist triangulated equivalences
		$$\underline{\GF_n(R)\cap \cot_n(R)} \backsimeq\cdots \backsimeq \underline{\GF_1(R)\cap \cot_1(R)} \backsimeq \underline{\GF(R)\cap \cot(R)}\simeq \underline{\PGF(R)}.$$
		Furthermore, if $\GPD_{AC}(R)<\infty$, then these triangulated categories are compactly generated.
		
	\end{cor}
	\proof  By Theorem \ref{GFn mod struc} and \cite[Theorem 4.9]{SS20} we have the model structures
	$$\M_n=(\GF_n(R),\PGF(R)^\perp,\cot_n(R)) \text{ and } \mathcal{N}=(\PGF(R), \PGF(R)^\perp, R\text{-Mod})$$
	with the same trivial objects and the following inclusions
	$$\PGF\subseteq \GF(R)\subseteq \GF_1(R)\subseteq \cdots \subseteq \GF_n(R).$$
	Applying \cite[Lemma 5.4]{EG19}, we get the following triangulated equivalences
	$$\Ho(\M_n) \backsimeq\cdots \backsimeq \Ho(\M_1) \backsimeq \Ho(\M_0)\simeq \Ho(\mathcal{N}).$$
	
	On the other hand, following Corollary \ref{quas-Frob} and \cite[page 27]{SS20}, we have 
	$$\Ho(\M_n)\backsimeq \underline{\GF_n(R)\cap \cot_n(R)} \text{ and }\Ho(\mathcal{N})\backsimeq \underline{\PGF(R)}.$$
	as triangulated categories. Hence, the first statement follows.

    Finally, if $\GPD_{AC}(R)<\infty$, then our last claim follows by \cite[Theorem 24, Proposition 21 and Corollary 33]{LW20} and the above triangulated equivalences.
	\cqfd
	
	It is well-known that the stable categories $\underline{\GP(R)}$ and $\underline{\GI(R)}$ are triangulated categories. Li and Wang showed in \cite[Theroem 1]{LW20} that these triangulated categories are compactly generated whenever $\GPD(R)<\infty$ with $R$ right or left coherent. Their proof is mainly based on \cite[Theorem 25]{LW20}. However,  the coherence assumption in \cite[Theorem 25]{LW20} is not needed if we replace $\GPD(R)<\infty$ with $\GPD_{AC}(R)<\infty$. Therefore, following the same argument in \cite[Corollary 34]{LW20}, we have the improved result:
	
	\begin{cor} Let $R$ be a ring with $\GPD_{AC}(R)<\infty$. Then,
		$$\underline{\GP(R)}\simeq \underline{\GI(R)}\simeq \underline{\GF(R)\cap\cot(R)}$$
		are compactly generated.
		
	\end{cor}

	\section{PGF dimension of modules and rings}
	
	\hskip .5cm In this section, we are interested in the PGF dimension of modules and rings and its connection with other homological dimensions.
	
	The following result is a key result to prove many of the results that will follow. It is a combination of \cite[Theorem 3.20]{Ste22} and \cite[Theorem 3.4$(i)\Leftrightarrow (ii)$]{DE22}, except for the last statement.
	\begin{prop}\label{charac of PGFn} The following assertions are equivalent:
		\begin{enumerate}
			\item $\PGFpd_R(M)\leq n.$
			\item $M$ is a direct summand in a module $N$ such that there exists an exact sequence of $R$-modules 
			$$0\to N\to X\to N\to 0$$ where $\pd_R(X)\leq n$ and $\Tor_{n+1}^R(E,N)=0$ for any injective right $R$-module.
			\item There exists an exact sequence of $R$-modules $$0\to M\to X\to G\to 0$$  where $G\in\PGF(R)$ and $\pd_R(X)\leq n$.
		\end{enumerate}
		
		Consequently, $\PGF_n(R)\cap \PGF(R)^{\perp}=\proj_n(R).$
	\end{prop}
	\proof (1) is equivalent to (2) by  \cite[Theorem 3.20]{Ste22} and  to (3) by \cite[Theorem 3.4]{DE22}. We only have to prove the equality $\PGF_n(R)\cap \PGF(R)^{\perp}=\proj_n(R)$.
	
	Clearly $\proj_n(R)\subseteq \PGF_n(R)$ and by Lemma \ref{n-flat is PGC}, $\proj_n(R)\subseteq \PGF(R)^\perp$. Then, $\proj_n(R)\subseteq \PGF_n(R)\cap \PGF(R)^{\perp}$. Conversely, if $M\in\PGF_n(R)\cap \PGF(R)^{\perp}$, then there exists by (3) a split exact sequence of $R$-modules 
	$0\to M\to X\to G\to 0$ where $G\in\PGF(R)$ and $\pd_R(X)\leq n$. Then, $M\in \proj_n(R)$.
	\cqfd

	As a consequence, we get new characterizations of left ($n$-)perfect rings. Recall that \cite[Definition 2.1]{EJL05} a ring $R$ is left $n$-perfect if every flat $R$-module has projective dimension at most  $n$. In particular, left perfect rings are exactly left $0$-perfect rings.
	
	\begin{thm}\label{n-perf}
		A ring $R$ is left $n$-perfect if and only every Gorenstein flat has PGF dimension at most $n$. In particular, the following assertions are equivalent:
		\begin{enumerate}
			\item $R$ is left perfect.
			\item $\GF(R)=\PGF(R)$.
			\item $\PGF(R)$ is closed under direct limits.
			\item $\left(\PGF(R),\PGF(R)^\perp\right) $ is a perfect cotorsion pair.
			\item $\PGF(R)$ is covering.
			\item Every Gorenstein flat $R$-module has  a PGF cover.
			\item Every flat $R$-module has a projective cover.	
		\end{enumerate}
		
	\end{thm}
	\proof $(\Leftarrow)$ Clearly $\flat(R)\subseteq \GF(R)\subseteq \PGF_n(R)$ and by Lemma \ref{n-flat is PGC}, $\flat(R)\subseteq \PGF(R)^\perp$. Thus, $\flat(R)\subseteq \PGF_n(R)\cap \PGF(R)^\perp=\proj_n(R)$ by Proposition \ref{charac of PGFn}. This means that $R$ is left $n$-perfect. $(\Rightarrow )$ Let $M$ be a Gorenstein flat $R$-module. By \cite[Theorem 3.5 and Proposition 3.6]{BM07}, $M$ is a direct summand in a module $N$ such that there exists an exact sequence of $R$-modules  $0\to N\to F\to N\to 0$
	with $F$ flat and $\Tor_{1}^R(E,N)=0$ for any injective right $R$-module. It follows that $\Tor_{k\geq 1}^R(E,N)\cong \Tor_{1}^R(E,N)=0$. Since $F\in\flat(R)\subseteq \proj_n(R)$, it follows from Proposition \ref{charac of PGFn} that $M\in \PGF_n(R)$.  
	
	We now prove the equivalences (1)-(7):
	
	$(1)\Leftrightarrow (2)$ By \cite[Theorem 4.9]{SS20}. It also follows from the above equivalence. 
	
	$(2)\Rightarrow (3)$ By Proposition \ref{Gn-flat prop}(3). 
	
	$(3)\Rightarrow (4)$ By  \cite[Theorem 4.9]{SS20} and \cite[Corollary 2.3.7]{GT12}. 
	
	$(4)\Rightarrow (5)$ Clear. 
	
	$(5)\Rightarrow (6)$ Clear. 
	
	$(6)\Rightarrow (7)$ Let $F$ be a flat $R$-module and $f:P\to F$ be a PGF cover. We claim that $f$ is a projective cover. By Wakamatsu Lemma, $K:=\Ker f\in \PGF(R)^\perp$. But $F\in\flat(R)\subseteq \PGF(R)^\perp$ by Lemma \ref{n-flat is PGC} and $\PGF(R)\cap \PGF(R)^\perp=\proj(R)$ by Proposition \ref{charac of PGFn}. Hence, $P\in \PGF(R)\cap \PGF(R)^\perp=\proj(R)$. It remains to show that for any morphism $g: P\to P$ such that $fg=f$, is an automorphism. But this holds true as $f$ is a PGF cover. 
	
	$(7)\Rightarrow (1)$ By the proof of \cite[Theorem 5.3.2(3)$\Rightarrow$ (1)]{EJ00}. \cqfd

	Theorem \ref{n-perf} can be interpreted from the point of view of global dimensions. For this purpose (and for later use), we introduce the following homological invariant.

	\begin{defn}The global PGF dimension of $R$ with respect to a class $\X\subseteq R$-Mod is defined as:
		$\PGFD_{\X}(R):=\sup\{\PGFpd_R(X)|X\in\X\}.$
		
		In particular, we set 
		$$\PGFD(R):=\PGFD_{R\text{-Mod}}(R)=\sup\{\PGFpd_R(M)| M\in R\text{-Mod}\}$$
		and we simply call it the global PGF dimension of $R$. 
	\end{defn}
	We note that the global PGF dimension has been recently  studied by Dalezios and Emmanouil \cite{DE22}.
	
	\begin{cor}A ring $R$ is left $n$-perfect if and only if  $\PGFD_{\GF(R)}(R)\leq n$.
		
		In particular, $R$ is left perfect if and only if  $\PGFD_{\GF(R)}(R)=0$.
	\end{cor}

	We now turn our attention to the questions raised in the Introduction. 
	\begin{rem} \label{GP1,2,3} Consider the following assertions:\\
		(GP1) Any Gorenstein projective $R$-module is Gorenstein flat.\\
		(GP2) Any $R$-module has a special Gorenstein projective precover.\\
		(GP3) Any Gorenstein projective $R$-module is PGF, that is, $\PGF(R)=\GP(R).$ 
		
		It follows that  (GP3) is equivalent to (GP1) by \cite[Theorem 3]{Iac20} and implies  (GP2) by \cite[Theorem 4.9]{SS20}. However, it is not clear whether (GP2) implies (GP1) or (GP3). Following this remark, our focus will be on (GP3).
	\end{rem}

	Instead of asking whether any Gorenstein projective module is PGF or not, one could ask the following question:
	
	(Q4): How far is a Gorenstein projective module from being PGF?
	
	\bigskip
	
	The following result reduces this problem to two cases. 
	\begin{prop} \label{PGF-pd(GP)} The PGF dimension of a Gorenstein projective $R$-module $M$ is either zero or infinite.
		
		
		Consequently, the global dimension $\PGFD_{\GP(R)}(R)$ is either zero or infinite.
		
	\end{prop}
	\proof Assume  $n=\PGFpd_R(M)<\infty$ and  let us consider a special PGF precover (which exists by \cite[Theorem 4.9]{SS20}) $0\to N\to G\to M\to 0.$
	By \cite[Proposition 2.4(ii)]{DE22}, $\PGFpd_R(N)\leq max\{\PGFpd_R(G),\PGFpd_R(M)\}\leq n.$
	Hence, $N\in\PGF_n(R)\cap\PGF(R)^\perp=\proj_n(R)$ by Proposition \ref{charac of PGFn}. Moreover, since $M\in\GP(R)$, there exists an exact sequence of $R$-modules $0\to M\to P_{-1}\to \cdots \to P_{-n}\to L\to 0$ with each $P_j$ projective. Hence, $\Ext^1_R(M,N)\cong \Ext^{n+1}_R(L,N)=0$ as $N\in\proj_n(R)$ and then the above short exact sequence splits. Therefore, $M\in \PGF(R)$ and this means that $\PGFpd_R(M)=0$.   \cqfd

	In order to continue our investigations, we need the notion of definable classes. Recall from \cite[Theorem 3.4.7]{Pre09} that a class of modules $\mathcal{D}$ is called definable if it is closed under direct products, direct limits and pure submodules. Such a class is, in particular, closed under direct sums and direct summands. We denote by $\langle\X\rangle$ the definable class generated by a class of modules $\X$, that is, the smallest definable class containing $\X$. The definable closure $\langle \X\rangle$ can be constructed, for instance, by closing $\X$ under direct products, then under pure submodules and, finally, under pure quotients. In case $\X=\{X\}$, we simply write $\langle\X\rangle=\langle X\rangle$.

	It has been shown in many recent works \cite{EIP20,CS21,SS20} that the definable classes $\langle R\rangle$ and $\langle R^+\rangle$ are of great importance. They will also be of interest to us. Thus, we propose the following definition, giving names to the modules belonging to these classes.

	\begin{defn}An $R$-module is called definable flat if it belongs to the definable class $\langle  R\rangle$. Dually, an  $R$-module is called definable injective if it belongs to the definable class $\langle  (R_R)^+\rangle$. 
	\end{defn} 
	
	Clearly,  $\langle R\rangle=\langle \proj(R)\rangle=\langle \flat(R)\rangle$ and $\langle (_RR)^+\rangle=\langle\inj(R^{op})\rangle =\langle\FPinj(R^{op})\rangle $. Moreover, $R$ is right coherent if and only if definable flat $R$-modules coincide with flats if and only if definable injective right $R$-modules coincide with FP-injectives (see \cite[Theorem 3.4.24]{Pre09}).
	
	As we will see next, the class of definable flat modules is a key point for extending many results from the class of coherent rings to a larger class of rings. 
	
	\begin{lem}\label{<R>=Fn(R)} Given an integer $n\geq 0$, $\fd_R(\langle R\rangle)\leq n$ if and only if the direct product of any family of flat $R$-modules has finite flat dimension $\leq n$. 
		
		Consequently, $\fd_R(\langle R\rangle)=\sup\{\fd_R(\prod_iF_i)|  \text{ $(F_i)_i$ is a family of flat $R$-modules} \}.$
	\end{lem}
	\proof ($\Rightarrow $) For every family $(F_i)_i$  of flat $R$-modules, we have $F_i\in\langle \flat(R)\rangle=\langle R\rangle$ and since $\langle R\rangle$ is closed under products, $\prod_i F_i\in \langle R\rangle$. Hence, $\fd_R(\prod_iF_i)\leq\fd_R(\langle R\rangle)\leq n$. ($\Leftarrow$) Assume $\fd_R(\prod\flat(R))\leq n$  and let $X\in \langle R\rangle=\langle \flat(R)\rangle$. The definable class $\langle \flat(R)\rangle$ can be constructed by closing $\flat(R)$ under products, then under pure submodules and finally under pure quotients.  By hypothesis, any direct product of any family of flat $R$-modules belongs to $\flat_n(R)$, and clearly $\flat_n(R)$ is closed under pure submodules and quotients (since a short exact sequence $\mathcal{E}$ is pure if and only if its character $\mathcal{E}^+$ is split). Hence, $\langle R\rangle=\langle \flat(R)\rangle \subseteq\flat_n(R)$, that is, $\fd_R(\langle R\rangle)\leq n$. \cqfd
	
	Izurdiaga \cite{CI16} investigated rings for which the direct product of any family of flat $R$-modules has finite flat dimension. He calls them \cite[Section 4]{CI16} right weak coherent rings. Moreover, if $n$ is the maximum of the set consisting of all flat dimensions of all direct products of flat $R$-modules,  then  $R$ is called right weak $n$-coherent. In particular,  right weak $0$-coherent rings are nothing but right coherent rings.  It follows by Lemma \ref{<R>=Fn(R)} that $R$ is right weak coherent if and only if and only if any definable flat module has finite flat dimension. In this case, $R$ is right weak $n$-coherent with $n=\fd_R(\langle R\rangle)$.

	We also recall from \cite{Gil10} the notion of Ding projective modules (which are introduced for the first time in \cite{DLM09} under a different name). An $R$-module $M$ is called Ding projective if it is a syzygy of a $\Hom_R(-,\flat(R))$-exact exact sequence of projective $R$-modules. The Ding projective dimension of an $R$-module $M$ is defined as ${\rm Dpd}_R(M)=\resdim_{\mathcal{DP}(R)}(M)$ where $\mathcal{DP}(R)$ denotes the class of Ding projective $R$-modules. 
	
	As a consequence of the above lemma, the assumption that $R$ is right coherent in \cite[Theorem 2]{Iac20} and \cite[Corollary 5.10]{CS21}  can be replaced by that $R$ is right weak coherent. 
	
	\begin{cor}\label{PGF=DP} Assume that $R$ is a right weak coherent ring. Given an $R$-module $M$, we have $\Dpd_R(M)= \PGFpd_R(M)$. In particular, $\DP(R)=\PGF(R)$ is the left hand of a complete cotorsion pair.
	\end{cor}
	\proof We only need to show the equality $\PGF(R)=\DP(R)$. The inclusion $\PGF(R)\subseteq \DP(R)$ holds for any ring $R$. For the other inequality, let $M\in\DP(R)$, then $M$ is a syzygy of a $\Hom_R(-,\flat(R))-$exact exact sequence $\PP$ of projective $R$-modules. Let $X$ be a definable flat $R$-module. By Lemma \ref{<R>=Fn(R)}, $m=\fd_R(X)<\infty$. Then, $X$ has a finite flat resolution $$0\to F_m\to \cdots \to F_0\to X\to 0.$$  
	
	Applying the functor $\Hom_R(\PP,-)$ to this exact sequence, we get an exact sequence of complexes
	$$0\to \Hom_R(\PP,F_m)\to \cdots \to\Hom_R(\PP,F_0) \to \Hom_R(\PP,X)\to 0.$$ As the class of exact complexes is thick, it follows that the complex $\Hom_R(\PP,X)$ is exact and hence $M\in\PGF(R)$.
	
	The last claim follows by \cite[Theorem 4.9]{SS20}. \cqfd

	One could, as in Lemma \ref{<R>=Fn(R)}, consider rings for which every definable flat has finite projective dimension. The following lemma, which we will use later, says that such rings are nothing new.
	
	\begin{lem}\label{pd(<R>)} $\pd_R(\langle R\rangle)<\infty$ if and only if $R$ is left weak coherent and right $n$-perfect for some integer $n\geq 0$. 
	\end{lem}
	\proof $(\Rightarrow)$ Assume $n=\pd_R(\langle R\rangle)<\infty$. Then, $\pd_R(X)\leq n$ for all $X\in \langle R\rangle$. In particular, $\fd_R(X)\leq \pd_R(X)\leq n$ for all $X\in \langle R\rangle$ and $\pd_R(X)\leq n$ for all $X\in\flat(R)\subseteq \langle R\rangle$.  
	$(\Leftarrow)$ Assume that $m=\fd_R(\langle R\rangle)<\infty$. Given a partial projective resolution of an $R$-module $M\in \langle R\rangle$
	$$0\to K_m\to P_{m-1}\to \cdots \to P_0\to M\to 0,$$ 
	we get that $K_m$ is flat and hence $\pd_R(K_m)\leq n$ as $K_m\in \langle R\rangle$. Therefore, $\pd_R(M)\leq n+m$. Thus, $\pd_R(\langle R\rangle)\leq n+m<\infty$. \cqfd

	Now we are ready to give a partial answer to Question (Q3)  raised in the Introduction. Recall that a $\Hom_R(-,\proj(R))$-exact exact complex of projective modules is called totally acyclic complex of projectives.

	\begin{prop}\label{PGF=GP}  $\PGF(R)=\GP(R)$ if and only if any totally acyclic complex of projectives is $\Hom_R(-,\langle R\rangle)$-exact if and only if any totally acyclic complex of projectives is $(\inj(R^{op})\otimes-)$-exact.

		In particular, this is the case if one of the following holds:  
		\begin{enumerate}
			\item[(a)] $R$ is right weak coherent and left $n$-perfect for some integer $n\geq 0$.
			
			\item[(b)] $R$ is a ring where any injective right $R$-module has finite flat dimension.
			
		\end{enumerate}
		
	\end{prop}
	\proof 	The equivalences follow immediately from \cite[Corollary 4.5]{SS20}.
	
	(a) By above, it suffices to show that any totally acyclic complex of projectives $\PP$ is $\Hom_R(-,X)$-exact for any definable flat $R$-module $X$. But, this can be shown as in the proof of corollary \ref{PGF=DP}, using a finite projective resolution of $X$ which exists by assumption and Lemma \ref{pd(<R>)}. 
	
	(b) Follows by \cite[Proposition 9]{Iac20}.
	\cqfd

	\begin{rem}\label{coh to weak coh}\item  
		\begin{enumerate}
			\item[(i)] Notice that the assertion "any totally acyclic complex of projectives is $\Hom_R(-,\langle R\rangle)$-exact" is equivalent to $\langle R\rangle\subseteq \GP(R)^\perp$. Therefore, Proposition \ref{PGF=GP} improves \cite[Theorem 4]{Iac20} in the sense that if we replace the class of flat modules with that of definable flat we can drop the coherence assumption. We also note that every right coherent ring is right weak coherent, so the example in Proposition \ref{PGF=GP}(a) is more general.
			
			\item[(ii)] In view of  Proposition \ref{PGF=GP}, Remark \ref{GP1,2,3} and Theorem \ref{n-perf}, if $R$ is right weak coherent, then  $R$ is left perfect if and only if $\GP(R)=\GF(R)$. This result improves \cite[Theorem 4.5 $(i)\Leftrightarrow (iv)$]{CET20}.
		\end{enumerate} 
	\end{rem}
	
	Consequently, partial answers to questions (Q1) and (Q2) are obtained.
	
	\begin{cor} \label{Q1 and Q2} Assume that $R$ satisfies the equivalent conditions of Proposition \ref{PGF=GP}.
		\begin{enumerate}
			\item For any $R$-module $M$, we have  $\Gfd_R(M)\leq \Gpd_R(M)$.
			
			In particular, every Gorenstein projective $R$-module is Gorenstein flat.
			\item  $(\GP(R),\GP(R)^\perp)$ is a complete hereditary cotorsion pair.
			
			In particular, every $R$-module has a special Gorenstein projective precover.
		\end{enumerate}
	\end{cor}
	\proof 1. We only need to show that $\GP(R)\subseteq \GF(R)$. But this follows by Proposition \ref{PGFD=GPD} as $\GP(R)=\PGF(R)\subseteq \GF(R)$.
	
	2. Follows by Proposition \ref{PGF=GP} and \cite[Theorem 4.9]{SS20}.
	\cqfd

	We now focus our attention on the global PGF dimension of $R$. As a first result, we provide simple ways to compute it (compare with \cite[Theorem 5.1]{DE22}).
	
	\begin{thm}\label{charc of PGFD}The following assertions are equivalent:
		\begin{enumerate}
			\item $\PGFD(R)\leq n$.
			\item  The following two assertions hold.
			\begin{enumerate}
				\item $\id_R(M)\leq n$ for every definable flat $R$-module $M$. 
				\item $\pd_R(M)\leq n$ for every injective $R$-module $M$.
			\end{enumerate}
			\item The following two assertions hold.
			\begin{enumerate}
				\item $\fd_{R^{op}}(M)\leq n$ for every definable injective right $R$-module $M$.
				\item $\pd_R(M)\leq n$ for every injective $R$-module $M$.
			\end{enumerate}
			\item The following two assertions hold.\begin{enumerate}
				\item $\fd_{R^{op}}(M)\leq n$ for every (FP-)injective right $R$-module $M$.
				\item $\pd_R(M)\leq n$ for every injective $R$-module $M$.
			\end{enumerate}
			
		\end{enumerate}
		
		Consequently, the global PGF dimension of $R$ can be computed by the following formulas:\begin{eqnarray*}
			\PGFD(R)&=&max\{\pd_R(\inj(R)),\id_R(\langle R\rangle)\}\\
			&=&max\{\pd_R(\inj(R)),\fd_{R^{op}}(\langle R^+\rangle)\}\\
			&=& max\{\pd_R(\inj(R),\fd_{R^{op}}(\inj(R^{op}))\}.
		\end{eqnarray*}
	\end{thm}
	
	\proof $1.\Rightarrow 2.$	(a) Consider an $R$-module $N$ and a finite PGF resolution:
	$$0\to G_n\to \cdots \to G_0\to N\to 0.$$
	
	By \cite[Corollary 4.5]{SS20}, $\Ext^{k\geq 1}_R(G_i,M)=0$ for any  $M\in\langle R\rangle$ and $i$. Hence, $\Ext^{n+1}_R(N,M)\cong \Ext_R^1(G_n,M)=0$. Then, $\id_R(M)\leq n$.
	
	(b) Let $M$ be an injective $R$-module. Since $\PGFpd_R(M)\leq n$, there exists a split short exact sequence $0\to M\to X\to G\to 0$ with $\pd_R(X)\leq n$ by Proposition \ref{charac of PGFn}(3). Hence, $\pd_R(M)\leq n$. 
	
	$2.\Rightarrow 3.$ We only prove (a) as (b) is clear. For any definable injective right $R$-module $M$, we have 	$\fd_{R^{op}}(M)\leq n$ if and only if  $\Tor^R_{n+1}(M,-)=0$ if and only if  $ \Ext_R^{n+1}(-,M^+)\cong\Tor^R_{n+1}(M,-)^+=0$ if and only $\id_R(M^+)\leq n$. This later holds by 2(a)  and \cite[Remark 2.10]{EIP20} as $M^+\in\langle R\rangle$. Hence, $\fd_{R^{op}}(M)\leq n$.

	$3.\Rightarrow 4.$ (b) is clear and (a) follows by $\inj({R^{op}})\subseteq \FPinj({R^{op}})\subseteq \langle \FPinj({R^{op}})\rangle= \langle (R_R)^+\rangle$.
	
	$4.\Rightarrow 1.$  Let $M$ be an $R$-module. Consider a projective and an injective resolution of $M$: $$\cdots \to P_1\to P_0\to M\to 0 \text{ and } 0\to M\to I_0\to I_1\to\cdots,$$ 
	respectively. Decomposing these exact sequences into short exact ones we get, for every integer $i\in \N$, 
	$$0 \to L_{i+1}\to P_i\to L_i\to 0\text{ \;\;and\;\;  }0\to K_i\to I_i\to K_{i+1}\to 0$$ 
	where $L_i=\Coker(P_{i+1}\to P_{i})$ and  $K_i=\Ker(I_i\to I_{i+1})$. Note that $M=L_0=K_0$. 
	Adding the direct sum of the first sequences, 
	$$0 \to \bigoplus_{i\in\N}L_{i+1}\to \bigoplus_{i\in\N}P_i\to M\oplus( \bigoplus_{i\in\N}L_{i+1})\to 0$$
	to the direct product of the second ones, 
	$$0\to  M\oplus( \prod_{i\in\N}K_{i+1})\to \prod_{i\in\N}I_i\to \prod_{i\in\N}K_{i+1}\to 0$$
	we get an exact sequence of the form
	$0\to N\to X\to N\to 0$, where $$X=(\bigoplus_{i\in N} P_i)\oplus (\prod_{i\in N} I_i)\text{ and } N=M\oplus \left((\bigoplus_{i\in N} L_{i+1}) \oplus (\prod_{i\in N} K_{i+1})\right).$$
	
	By (b), $\pd_R(X)=\pd_R((\prod_{i\in N} I_i))\leq n$ and by (a), $\Tor^R_{n+1}(E,N)=0$ for any injective right $R$-module $E$. And since $M$ is a direct summand in $N$, it follows from Proposition \ref{charac of PGFn} that $\PGFpd_R(M)\leq n$ as desired.
	\cqfd
	
	Global dimensions of rings are known to characterize classical rings. And the global PGF dimension is no exception. 
	Recall that a ring $R$ is called quasi-Frobenius if injective left (resp., right) $R$-modules coincide with projective left (resp., right) modules.
	
	\begin{cor}$\PGFD(R)=0$ if and only if $R$ is quasi-Frobenius.	
	\end{cor}
	\proof  $(\Rightarrow)$ By Theorem \ref{charc of PGFD}(2), we have $\proj(R)\subseteq \langle R\rangle\subseteq \inj(R)$ and $\inj(R)\subseteq \proj(R)$. Then, $\proj(R)=\inj(R)$ and therefore $R$ is quasi-Frobenius.  $(\Leftarrow)$ If $R$ is quasi-Frobenius, then $\inj(R)=\proj(R)$ and $\inj(R^{op})=\proj(R^{op})\subseteq \flat(R^{op})$. Hence, $\PGFD(R)=0$ by Theorem \ref{charc of PGFD}(4).
	\cqfd

	At this point, it is natural to ask what the relationship between the global PGF dimension and the global Gorenstein projective dimension could be. Clearly, $\GPD(R)\leq \PGFD(R)$ as $\PGF(R)\subseteq \GP(R)$. Over an (Iwanaga-)Gorenstein ring $R$, it is an immediate consequence of Huang \cite[Theorems 4.9 and 4.13]{Hua22} that we have an equality $\GPD(R)=\PGFD(R)$.  Dalezios and Emmanouil, on the other hand, have shown in \cite[Theorem 5.1]{DE22} this equality when $R$ has finite global PGF dimension. 
	
	Our last main result gives a positive answer to this question for any ring $R$. It is worth noting that the following result, as pointed out to me by  Sergio Estrada, is contained in the proof of a recent result by Wang, Yang, Shao, and Zhang \cite[Theorem 3.7]{WYSZ23}. However, our proof is different.

	\begin{thm} \label{PGFD=GPD}The global Gorenstein projective dimension of $R$ coincides with the global PGF dimension of $R$. That is,  $\GPD(R)=\PGFD(R).$
	\end{thm}
	\proof The inequality $\GPD(R)\leq \PGFD(R)$ is clear since $\PGF(R)\subseteq \GP(R)$. If $\GPD(R)=\infty$, then $\PGFD(R)=\GPD(R)$ and we are done.
	
	Assume $n=\GPD(R)<\infty$. By \cite[Corollary 2.7]{BM10}, $\pd_R(E)\leq n$ for every injective $R$-module $E$. On the other hand, following \cite[Examples 4.4(3)]{CI16}, $R$ is right weak coherent. It follows that $\fd_R(\langle R\rangle)<\infty$ by Lemma \ref{<R>=Fn(R)}. This implies that $\id_R(\langle R\rangle)\leq n$ by \cite[Corollary 2.7]{BM10}. Applying now Theorem \ref{charc of PGFD}(2), we get that $\PGFD(R)\leq n$.
	\cqfd
	We end the paper with some consequences of Theorem \ref{PGFD=GPD}.%
	
	Gorenstein injective modules are seen as dual to Gorenstein projective modules. Based on this, one could introduce a dual notion of a PGF module. Let us call an $R$-module Gorenstein definable injective if it is a syzygy of a $\Hom_R(\langle R^+\rangle,-)$-exact exact sequence of injective $R$-modules. Define its Gorenstein definable injective dimension as ${\rm Gdid}(M)=\resdim_{\mathcal{GDI}(R)}(M)$ where $\mathcal{GDI}(R)$ denotes the class of all Gorenstein definable injective $R$-modules. 
	
	It is known that the global dimension of $R$ can be computed either by projective or injective dimensions, that is, the following equality holds:
	$$\sup\{\pd_R(M)|\text{ \text{ $M\in R$-Mod}}\}=\sup\{\id_R(M)|\text{ $M\in R$-Mod}\}.$$
	
	This fact was extended to the Gorenstein setting by Enochs and Jenda \cite[Section 12.3]{EJ00} for (Iwanaga-)Gorenstein rings and later by Bennis and Mahdou \cite[Theorem 1.1]{BM10} for any ring. 
	
	Inspired by this, one could ask whether the equality

	$$\sup\{\PGFpd_R(M)|\text{$M\in R$-Mod}\}=\sup\{{\rm Gdid}_R(M)| \text{ $M\in R$-Mod}\}$$
	holds true as well? 
	
	Using again Theorem \ref{PGFD=GPD}, we give a positive answer to this question. 
	
	\begin{cor}\label{PGFD=GID=DID} For any ring $R$, we have the following equalities:
		$$\PGFD(R)=\GID(R)=\GDID(R)$$
		with $\GDID(R):=\sup\{\rm {Gdid}_R(M)|\text{ $M$ is an $R$-module}\}.$

	\end{cor}
	\proof  First, we have the equality $\PGFD(R)=\GID(R)$ that follows by Theorem \ref{PGFD=GPD} and \cite[Theorem 1.1]{BM10}.
	
	Clearly, $\GID(R)\leq \sup\{{\rm Gdid}(M)|\text{ $M\in R$-Mod}\}$ since $\mathcal{GDI}(R)\subseteq \GI(R)$.  Assume now that $n=\GID(R)<\infty$. If we prove that $\id_R(M)<\infty$ for any definable injective $R$-module $M$, we are done, as this condition implies that $\GI(R)=\mathcal{GDI}(R)$. Following \cite[Corollary 2.7]{BM10}, it suffices to show that $\fd_R(M)<\infty$ for any definable injective $R$-module. But, using \cite[Lemma 5.6(1)]{CS21}, it suffices to show that $\fd_R(M)<\infty$ for every FP-injective $R$-module $M$. 
	
	Let $_RM$ be FP-injective. Then, there exists a pure monomorphism $0\to M \hookrightarrow E$ with $_RE$ injective. Using again \cite[Corollary 2.7]{BM10}, we get that $\fd_R(E)\leq n$ and hence $\fd_R(M)\leq \fd_R(E)\leq n<\infty$ as desired.\cqfd

    \begin{rem}It clear by the definitions that we have the inequalities
    $$\GPD(R) \leq \DPD(R) \leq\PGFD(R)\leq \GPD_{AC}(R),$$
    $$\GID(R) \leq \DID(R) \leq\GDID(R)\leq \GID_{AC}(R)$$
      
with  $$\DPD(R):=\sup\{\Dpd_R(M)| \text{ $M$ is an $R$-module} \},$$
      $$\DID(R) :=\sup\{ \Did_R(M)| \text{ $M$ is an $R$-module} \},$$
      $$\GPD_{AC}(R):=\sup\{ \resdim_{\GI_{AC}(R)}(M)| \text{ $M$ is an $R$-module} \},$$
      where $\GI_{AC}(R)$ denotes the class of all Gorenstein AC-injective modules in the sense of \cite{BGH14}.

    As stated above, Bennis and Mahdou showed in \cite[Theorem 1.1]{BM10} the equality $\GPD(R)=\GID(R)$ for any ring $R$. Generally, Huerta, Mendoza and P\'erez have recently shown in \cite[Corollary 8.8]{HMP23} the following equalities:
    $$\GPD_{AC}(R)= \DPD(R)=\GPD(R)=\GID(R) =\DID(R)=\GID_{AC}(R).$$
    The result \cite[Corollary 8.8]{HMP23} is based on \cite[Lemma 8.1]{HMP23}, which is stated incorrectly as indicated in the corrigendum \cite{HMP24}. So, the result \cite[Corollary 8.8]{HMP23} is no longer available and the question of whether the equalities
$\DPD(R)=\DID(R)$ and $\GPD_{AC}(R)=\GID_{AC}(R)$ hold for any ring $R$ is still open. However, following Corollary \ref{PGFD=GID=DID}, we have a positive answer for the second equality. Generally, we have the following equalities for any ring $R$:
$$\PGFD(R)=\DPD(R)=\GPD(R)=\GID(R)=\DID(R)=\GDID(R).$$
    \end{rem}
\bigskip
 
 Auslander’s Theorem on the global dimension states that we can compute the global dimension of $R$ by just computing the projective dimension of cyclic $R$-modules. That is, the formula $\gldim(R)=\sup\{\pd_R(R/I)|\text{ $I$ is a left idea}\}$ holds true.  Bennis, Hu and Wang \cite[Theorem 1.1]{BHW15} extended this formula to the Gorenstein setting for commutative rings. However, one can see that the same proof also holds for non-commutative rings. Taking advantage of this fact, together with Theorem \ref{PGFD=GPD}, we get a PGF version of Auslander’s Theorem. 
	
	\begin{cor}\label{Aus}(\textbf{Auslander's Theorem on the global PGF  dimension}) For any ring $R$, we have the following formula:
		$$\PGFD(R)=\sup\{\PGFpd_R(R/I)|\text{ $I$ is a left ideal}\}.$$
	\end{cor}
	\proof The inequality $\{\PGFpd_R(R/I)|\text{ $I$ is a left ideal}\}\leq \PGFD(R)$ is clear. We may assume that   $n=\sup\{\PGFpd_R(R/I)|\text{ $I$ is a left ideal}\}<\infty$. Since $\PGF(R)\subseteq \GP(R)$, we have $\Gpd(R/I)\leq \PGFpd_R(R/I)\leq n$ for every left ideal $I$, and therefore $\sup\{\Gpd_R(R/I)|\text{ $I$ is a left ideal}\}\leq n$. On the other hand, by the non-commutative version of \cite[Theorem 1.1]{BHW15}, we have that $\GPD(R)=\sup\{\Gpd_R(R/I)|\text{ $I$ is a left ideal}\}\leq n$. Finally, using Theorem \ref{PGFD=GPD}, we get that $\PGFD(R)=\GPD(R)\leq n$ which completes the proof. \cqfd
	
	\bigskip
	\noindent\textbf{Acknowledgement.} 
 The author would like to thank Marco A. P\'erez for his kindness in sharing his book \cite{Per16} and for informing me about the Corrigendum \cite{HMP24}. Additionally, gratitude is extended to Sergio Estrada for identifying a misprint in the proof of Lemma 4.8 and for sharing the reference \cite{WYSZ23}. Special thanks also go to Luis Oyonarte for his comments on an earlier version of this paper and for his constant support. Finally, we express our appreciation to the referee for suggestions that enhanced the exposition.


\begin{thebibliography}{999}
	\small	

		\bibitem{AR94}
		 Ad\'{a}mek, J.,  Rosick\'{y}, J.: Locally presentable and accessible categories. Number 189 in London Mathematical Society Lecture Note Series, Cambridge University Press. (1994)
		
		\bibitem{AB69}
		 Auslander, M., Bridge,  M.: Stable Module Theory. Mem. Amer. Math. Soc. 94, Providence, RI.: Amer. Math. Soc. (1969)
		
		
		\bibitem{BEGO22} 
		 Bennis, D.,  El Maaouy, R.,  Garc\'{\i}a Rozas J.R.,  Oyonarte, L.: \textit{Relative Gorenstein flat modules and dimension}. Comm. Alg. \textbf{50}, 3853-3882 (2022)
		
		\bibitem{BEGO22b} 
		Bennis, D.,  El Maaouy, R.,  Garc\'{\i}a Rozas J.R.,  Oyonarte, L.: \textit{Relative Gorenstein flat modules and Foxby classes and their model structures}. arXiv:2205.02032 (2022)
		
		\bibitem{BHW15}
		 Bennis, D.,  Hu K.,  Wang, F.: \textit{Gorenstein analogue of Auslander’s theorem on the global dimension}. Comm. Alg.  \textbf{43}, 174-181 (2015)
		
		\bibitem{BM10}
		 Bennis, D.,  Mahdou, N.: \textit{Global Gorenstein dimensions}. Proc. Am. Math. Soc. \textbf{138},  461-465 (2010)
		
		\bibitem{BM07} 
		Bennis, D.,  Mahdou, N.: \textit{Strongly  Gorenstein  projective,  injective,  and  flat  modules}. J. Pure Appl. Algebra. \textbf{210}, 437-445 (2007)
		
		\bibitem{BGH14}
		 Bravo, D.,  Gillespie, J., Hovey, M.: \textit{The stable module category of a general ring}. https://arxiv.org/abs/1405.5768 (2014)
		
		
		\bibitem{CI16}
		 Cort\'{e}s-Izurdiaga, M.: \textit{Products of flat modules and global dimension relative to $\mathcal{F}$-Mittag-Leffler modules}. Proc. Amer. Math. Soc. \textbf{144}, 4557–4571 (2016)
		
		\bibitem{CS21}
		 Cort\'{e}s-Izurdiaga, M., \v{S}aroch, J.: \textit{Module classes induced by complexes and $\lambda$-pure-injective modules}. arXiv:2104.08602 (2021) 
		
		\bibitem{CET20}
		 Christensen, L. W., Estrada,  S.,  Thompson, P.: \textit{Homotopy categories of totally acyclic complex whit applications to the flat-cotorsion theory}. Contemporary Mathematics. \textbf{751}, 99–118 (2020)
		
		\bibitem{DE22}
		 Dalezios, G., Emmanouil, I.: \textit{Homological dimension based on a class of Gorenstein flat modules}. arXiv:2208.05692 (2022)
		
		
		\bibitem{DLM09}
		 Ding, N.,  Li, Y.,   Mao, L.: \textit{Strongly Gorenstein flat modules}. J. Aust. Math. Soc. \textbf{66}, 323–338  (2009)
		
		\bibitem{DM07}
		 Ding, N.,  Mao, L.: \textit{Envelopes and covers by modules of finite FP-injective and flat dimensions}. Comm. Alg. \textbf{35},  833–849  (2007)

		
		\bibitem{Emm12}
		  Emmanouil, I.:\textit{On the finiteness of Gorenstein homological dimensions}. J. Algebra. \textbf{372}, 376–396 (2012)
		
		\bibitem{EJ00}
		 Enochs, E.,  Jenda,  O.: Relative Homological Algebra. de Gruyter Expositions in Mathematics. \textbf{30}, Walter de Gruyter. (2000)
		
		
		\bibitem{EJL05}
		 Enochs, E.,  Jenda, O., L\'{o}pez-Ramos, J. A.: \textit{Dualizing modules and n-perfect rings}. Proc. Edinburgh Math. Soc. \textbf{48}, 75-90  (2005)
		
		\bibitem{EG19}
		 Estrada, S.,  Gillespie, J.: \textit{The projective stable category of a coherent scheme}. Proc. R. Soc. Edinb., Sect. A, Math. \textbf{149}, 15-43  (2019)
		
	
		\bibitem{EIP20} 
		 Estrada, S.,  Iacob, A., P\'{e}rez,  M. A.: \textit{Model structures and relative Gorenstein flat modules and chain complexes}. In Categorical, homological and combinatorial methods in algebra. \textbf{751}, 135-175  (2020)
		
		\bibitem{Gil11}
		 Gillespie, J.: \textit{Model structures on exact categories}. J. Pure. Appl. Algebra. \textbf{215}, 2892-2902 (2011)
		
		\bibitem{Gil16}
		 Gillespie, J.: \textit{Hereditary abelian model categories}. Bull. Lond. Math. Soc. \textbf{48}, 895-922  (2016)
		
		\bibitem{Gil10}
		 Gillespie, J.: \textit{Model structures on modules over Ding-Chen rings}. Homology Homotopy Appl, \textbf{12}, 61-73  (2010)
		
		
		\bibitem{GT12}
		 G\"{o}bel, R.,  Trlifaj,  J.:  \textit{Approximations and Endomorphism Algebras of Modules}. de Gruyter Expositions in Mathematics, \textbf{41}, 2nd revised and extended edition. Berlin Boston. (2012)
		
		
		\bibitem{Hov99}
		 Hovey, M.: \textit{Model categories}, Mathematical Surveys and Monographs 63. American Mathematical Society, Providence, RI. (1999)
		
		\bibitem{Hov02}
		 Hovey, M.: \textit{Cotorsion pairs, model category structures, and representation theory}. Math. Z. \textbf{241}, 553-592  (2002)
		
		\bibitem{Hua22} 
		  Huang, Z.: \textit{Homological dimensions relative to preresolving subcategories II}. Forum Math. \textbf{34}, 507–530  (2022)

         \bibitem{HMP24}
        Huerta, M. Y.,  Mendoza, O., P\'{e}rez, M. A.:  \textit{Corrigendum to “m-Periodic Gorenstein objects” [J. Algebra 621 (2023) 1–40]}, J. Algebra, 654, 70-81 (2024).
		\bibitem{HMP23}
		 Huerta, M. Y.,  Mendoza, O., P\'{e}rez, M. A.: \textit{$m$-periodic Gorenstein objects}. J. Algebra. \textbf{621}, 1-40 (2023)
        
		
		\bibitem{Iac20}
		 Iacob, A.: \textit{Projectively coresolved Gorenstein flat and ding projective modules}. Comm. Alg. \textbf{48}, 2883–2893  (2020)
		
		\bibitem{LW20}
		 Liang, L.,  Wang, J.: \textit{Relative global dimensions and stable homotopy categories}. C. R. Math. Acad. Sci. Paris. \textbf{358}, 379-392  (2020)
		
		\bibitem{Per16} 
		 P\'erez, M. A.: Introduction to Abelian Model Structures and Gorenstein Homological Dimensions. Monographs and Research Notes in Mathematics. CRC Press, Boca Raton, FL. (2016)
		
		\bibitem{Pre09}
	    Prest, 	M.: Purity, \textit{Spectra and Localisation}. \textbf{121}, of Encyclopedia of Mathematics and its Applications. Cambridge University Press, Cambridge, (2009)
		
		\bibitem{Ste22}
		 Stergiopoulou, D. D.:  \textit{Strongly n-projectively coresolved Gorenstein flat modules}. arXiv:2210.16816 (2022)
		
		\bibitem{SS20}
		 \v{S}aroch, J., \v{S}\v{t}ov\'{\i}\v{c}ek, J.: \textit{Singular compactness and definability for $\Sigma$-cotorsion and Gorenstein modules}. Selecta Math.  \textbf{26}, 23-40  (2020)

   \bibitem{WYSZ23}
   Wang, J.,  Yang, G.,  Shao, Q.,  Zhang, Q.:
   \textit{On Gorenstein global and Gorenstein weak global dimensions}. Colloq. Math. 174, 45-67 (2023)
		
		\bibitem{Xu17}
		 Xu., A. \textit{Gorenstein Modules and Gorenstein Model Structures}. Glasg. Math. J. \textbf{59}, 1-19  (2017)
		
		\bibitem{YL12}
		 Yang, G.,  Liu, Z.:  \textit{Gorenstein flat covers over GF-closed rings}. Commun. Alg. \textbf{40}, 1632–1640  (2012)
		
	\end{thebibliography}
\end{document}